\documentclass[a4paper,10 pt, conference]{article}  

\usepackage{graphicx}
\usepackage{cancel}
\usepackage{psfrag}
\usepackage{amsmath} 
\usepackage{url}

\usepackage{amssymb} 
\usepackage{color}

\title{{\bf A note about high-order semi-implicit differentiation: application to a numerical integration scheme with Taylor-based compensated error}}
\author{Loïc MICHEL$^{1}$ and Jean-Pierre BARBOT$^{1, \, 2}$}
\date{
\small{1. Nantes Université - Ecole Centrale de Nantes \\
Laboratoire des Sciences du Numérique de Nantes – UMR 6004 CNRS \\
1 rue de la Noë} \\
44300 Nantes, France \\
(corresponding author : \texttt{loic.michel54@gmail.com}) \\ 
2. QUARTZ EA 7393, ENSEA, Cergy-Pontoise, France}

\begin{document}
\maketitle

\vspace{2cm}
{\bf Abstract:} In this brief, we discuss the implementation of a  third order semi-implicit differentiator as a complement of the recent work by the author that proposes an interconnected semi-implicit Euler double differentiators algorithm through Taylor expansion refinement. The proposed algorithm is dual to the interconnected approach since it offers alternative flexibility to be tuned and to be implemented in real-time processes. In particular, an application to a numerical integration scheme is presented as the Taylor refinement can be of interest to improve the global convergence. Numerical results are presented to support the rightness of the proposed method.

\vspace{6cm}
Preprint - Work in progress

\newpage
\section{Introduction}

Acary {\it et al.} \cite{Acary:2012} introduced an implicit discretization technique from which differentiation techniques, for implicit differential inclusion, have emerged recently (see {\it e.g.} a survey in \cite{Rasool_CEP}). The principle of this method was initially to replace the classical sign function by an implicit projector, then it has been combined with explicit (homogeneous) terms to define the semi-implicit differentiation \cite{ifacWC}; the overall goal is to reduce the effects of chattering. 
Based on the latest work from the author {\it et al.}  \cite{VSS_2024_michel} that proposes a semi-implicit cascade differentiator including Taylor refinement to reduce the estimation error, in this brief, an alternative version of this latter algorithm is proposed that could possibly offer more flexibility towards practical implementation.

{The brief is organized as followed. The first Section defines the third order semi-implicit differentiator as a natural extension of the previous works \cite{VSS_2022_michel} \cite{VSS_2024_michel} including numerical results, the second Section illustrates an application towards a numerical integration scheme and the third Section gives some concluding remarks and perspectives.}
{Simulation code can be found in the repository: \texttt{https://github.com/LoicMichelControl/SemiImplicitDifferentiation.git}

\section{Semi-implicit differentiation algorithm including Taylor refinement}

\subsection{Presentation}

Considering a real signal $x_1(t)$ to differentiate, denote $x_1(k \delta)$ its corresponding discretization where $k$ is an integer and $\delta$ the sampling period. Moreover, $x_1^+$ stands for $x_1((k+1) \delta)$.
From the cascaded differentiator  \cite{VSS_2024_michel}, we derive a {\it third-order Semi-Implicit Differentiator} (SIHD-3), including {\it correction terms} (in bold) that refine the Taylor expansion of the two first estimation lines:

\begin{eqnarray} 
\left\{\begin{array}{l}
{z}_4^+= z_4+  E_1^+ E_2^+ E_3^+ h \left( \lambda_4 { \mu^4} \left|{e_1}\right|^{4\,\alpha-3} \mathcal{N}_4  \right)  \\
{z}_3^+= z_3+ E_1^+ E_2^+ h \left(z_4^++\lambda_3 { \mu^3} \left|{e_1}\right|^{3\,\alpha-2} \mathcal{N}_3 \right)  \\
{z_2^+}= z_2+ E_1^+ h \left(z_3^+ \boldsymbol{-E_3^+ h \frac{1}{2} z_4^+ } + \lambda_2 { \mu^2} \left|{e_1}\right|^{2\,\alpha-1} \mathcal{N}_2  \right) \\
z_1^+= z_1+ h \left( z_2^+ \boldsymbol{-E_2^+  \frac{h}{2} z_3^+ + E_3^+ E_2^+ \frac{h^2}{3!} z_4^+ } +  \lambda_1 {\mu} |{e_1}|^{\alpha} \mathcal{N}_1 \right)  \\
\end{array} \right. \label{eq:semi-implicit-diff-cor} 
 \end{eqnarray}
 
 \noindent
 or, equivalently, introducing the SIHD related operator:
 
 \begin{equation}
 z_i = \frak{D}_f ^{(i)}(z, z^+, x_1), \quad i = 1...4
  \label{eq:semi-implicit-D} 
 \end{equation}
 
 \noindent
 where $z_i$ is the estimation associated to the $i$th line of the differentiator \eqref{eq:semi-implicit-diff-cor}.
 
 \noindent
Denote $y = x_1$ the signal to differentiate; $\dot{y} = {x_2}$,  $\ddot{y} = x_3$ and $y^{(4)} = x_4$ are respectively the first, second and third order differentiation;  $ z_2, \, z_3$ and $z_4$ are respectively the first, second and third order related estimations; $e_1 = x_1 - z_1$ is the estimation error.

\newpage
\noindent
The projector\footnote{The projector behaves either as a sign function, or as a nonlinear "slope" depending on the value of the error $e_1$ \cite{ifacWC}.} $\mathcal{N}_q$ for $q = 1...4$  is defined by:

 \begin{eqnarray}
\mathcal{N}_q(\epsilon_1):= \left\{\begin{array}{l}
\epsilon_1 \in SD \hspace{2mm}   \rightarrow  \mathcal{N}_q={\dfrac{\lceil \epsilon_1\rfloor^{q(1-\alpha)}}{ \lambda_q (\mu h)^q} }, \hspace{2mm} E_{q}^{+}=1\\ \\
\epsilon_1 \notin SD \hspace{2mm}  \rightarrow  \mathcal{N}_q=   \mathrm{sign}(\epsilon_1), \hspace{2mm}  E_{q}^{+}=0
\end{array} \right. \label{projectorN1}
\end{eqnarray}
with the convergence domain\footnote{Note that, according to \cite{CDC_2021_michel}, it has been shown that the error $e_2$, corresponding to the second estimation line of \eqref{eq:semi-implicit-diff-cor}, converges asymptotically to $e_ 2^+ = e_1^+ / h$, hence 'propagating' the errors of estimation through the differentiation orders by 'increment' of $h$.} $SD=\{\epsilon_1 \, / \, |\epsilon_1|\leq (\lambda_1 \mu h)^{\frac{1}{q(1-\alpha)}}\}$.

In eq. \eqref{eq:semi-implicit-diff-cor}, the variable $E^+_q$ allows activating the upper $(q+1)$th line if the convergence domain has been reached at the $q$th line. This prevent from getting unexpected behavior between lines like oscillations that could destabilize the differentiator.

The corrections terms have been introduced by the author {\it et al.}  in \cite{VSS_2022_michel} with respect to the second order differentiation and a study of the convergence has been proposed in \cite{CDC_2021_michel}. The complete derivation of the refined Taylor expansion in \eqref{eq:semi-implicit-diff-cor}  is described in Appendix.

\subsection{Numerical results}

The  differentiator \eqref{eq:semi-implicit-diff-cor} is simulated over 10 sec. considering a sine function signal 

$$y(t) = \sin(2 \pi t )$$ 

\noindent
under the following initial conditions (that match the derivatives): $x_1(0) = 0, \, x_2(0) = 1, \, x_3(0) = 0$ and $x_4(0) = -1$.  The gains and parameters  are set as follows: $\lambda_1 = 10^3$, $\lambda_2 = 10^6$, $\lambda_3 = 10^9$, $\lambda_4 = 10^9$ and $\alpha = 0.95$.

\vspace{0.5cm}
\noindent
Figures \ref{fig_1} and \ref{fig_2} illustrate respectively the evolution of the estimated $z_1, z_2, z_3$ and $z_4$ with respect to the corresponding exact differentiation orders of the sine function without extra noise on $y$.

\noindent
Figures \ref{fig_3} and \ref{fig_4} illustrate respectively the evolution of the estimated $z_1, z_2, z_3$ and $z_4$ with respect to the corresponding exact differentiation orders of the sine function including extra noise of very small amplitude on $y$ (obtained by a random generator\footnote{The noise is computed under Matlab\textsuperscript{\tiny\textregistered} using the instruction \texttt{eta0 * rand()} where \texttt{eta0} is a factor that multiplies the noise.} that is multiplied by $\eta_0 = 10^{-6}$).

\begin{figure}[!b]
\begin{center}
   \includegraphics[width=13cm]{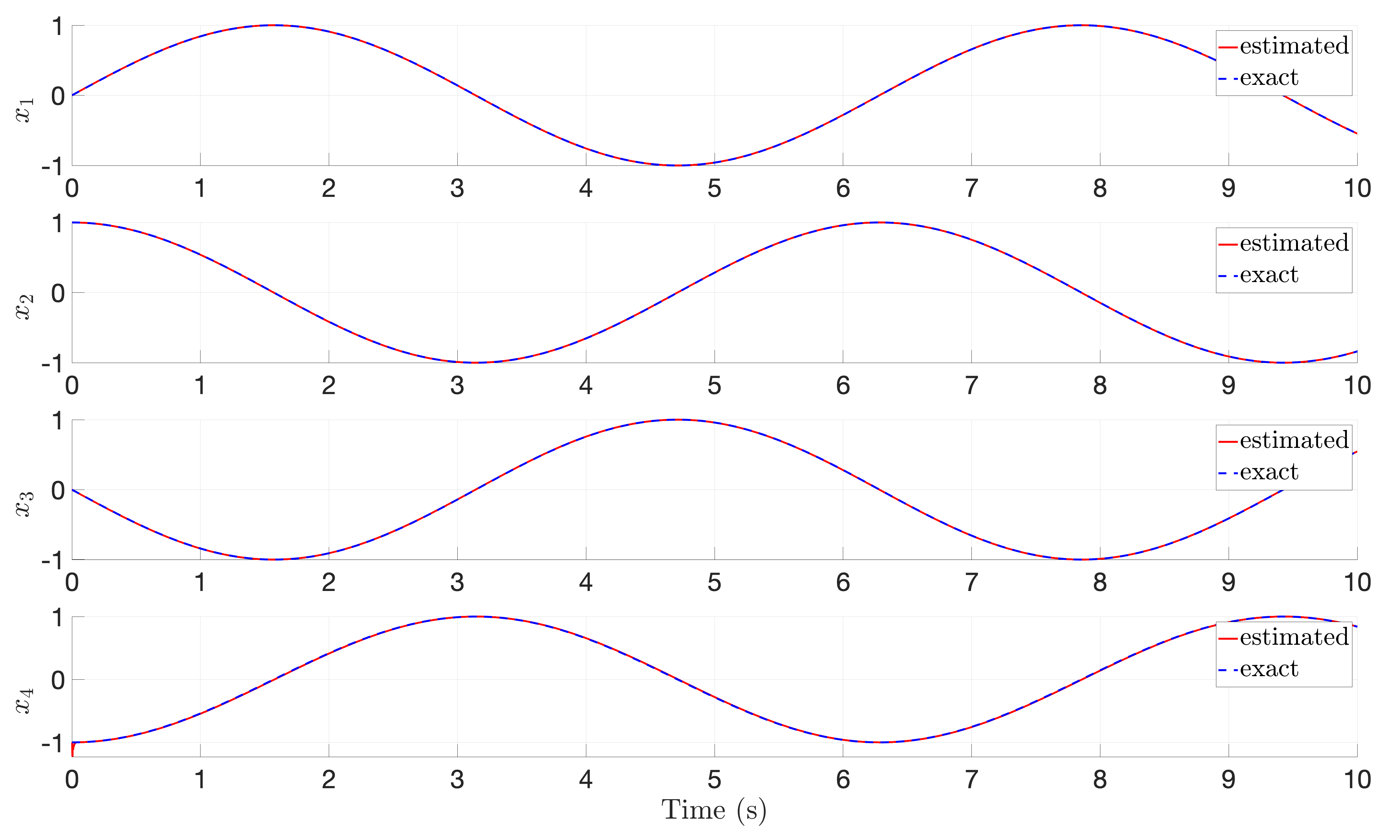}
      \end{center}
            \caption{Example of SIHD3-based differentiation of a sine function without extra noise.}
      \label{fig_1}
\end{figure}

\begin{figure}[!]
\begin{center}
   \includegraphics[width=13cm]{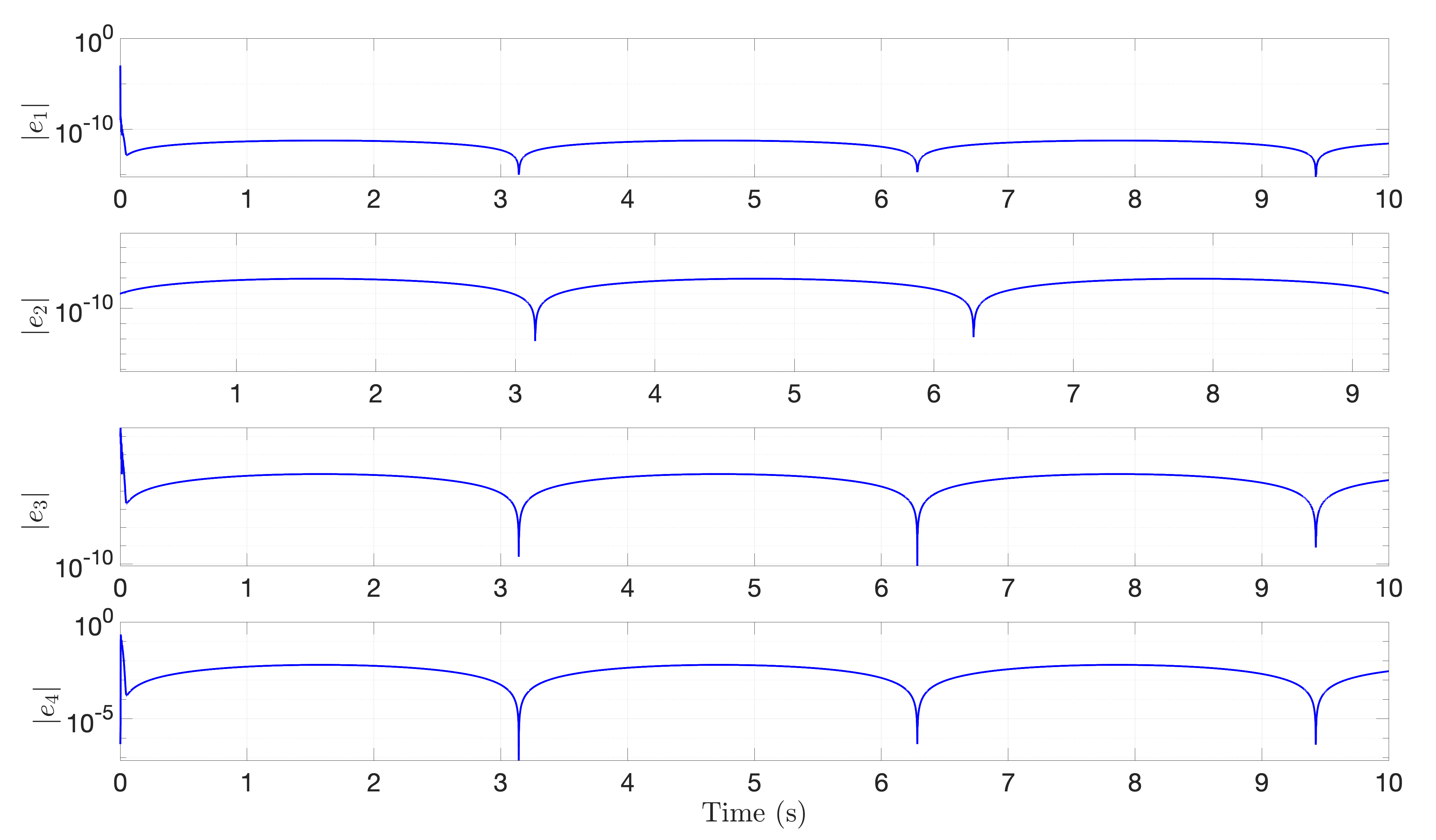}
      \end{center}
            \caption{Example of SIHD3-based differentiation of a sine function without extra noise - estimation errors.}

       \label{fig_2}
\end{figure}

\begin{figure}[!]
\begin{center}
   \includegraphics[width=13cm]{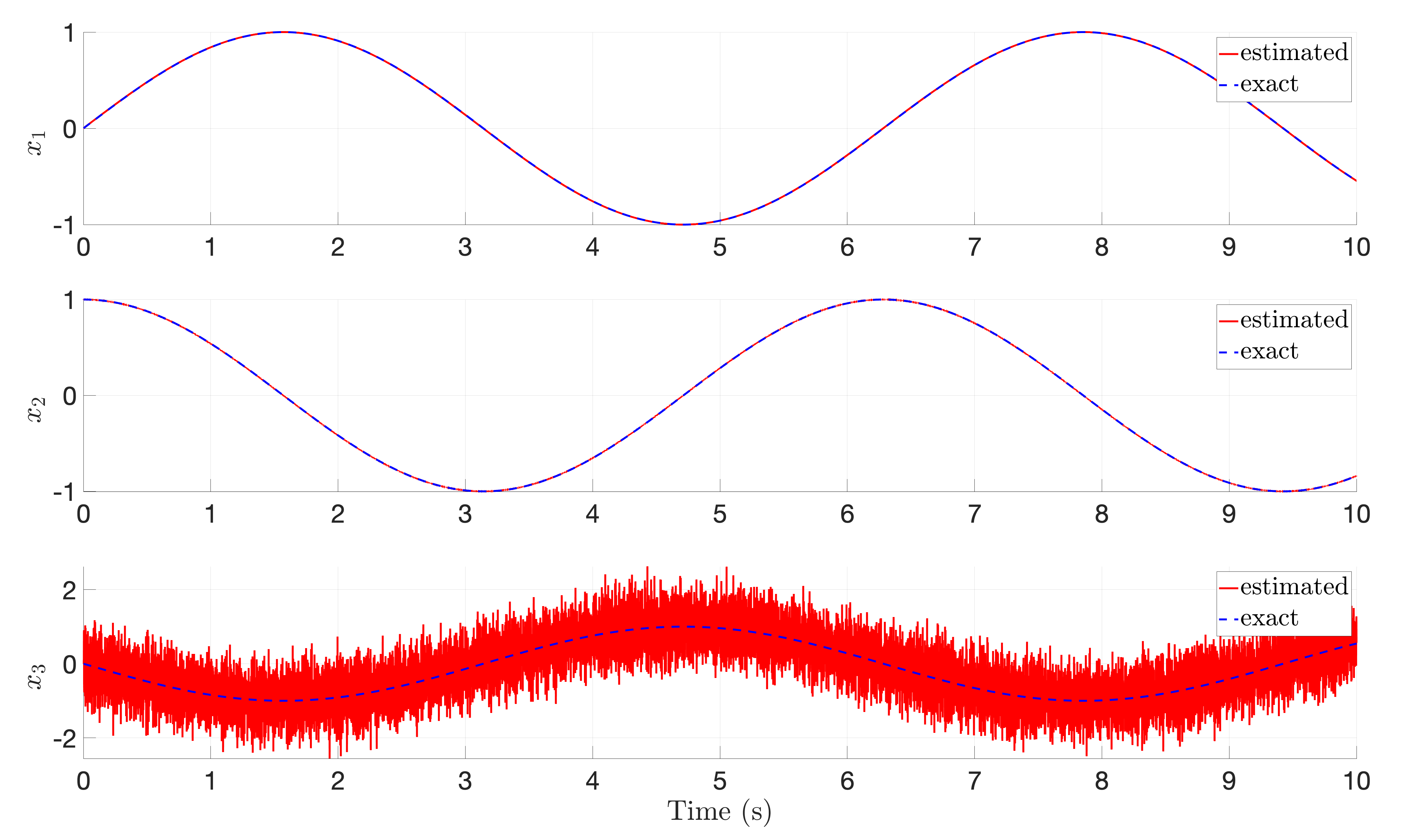}
      \end{center}
            \caption{Example of SIHD3-based differentiation of a sine function including extra noise.}
      \label{fig_3}
\end{figure}

\begin{figure}[!]
\begin{center}
   \includegraphics[width=13cm]{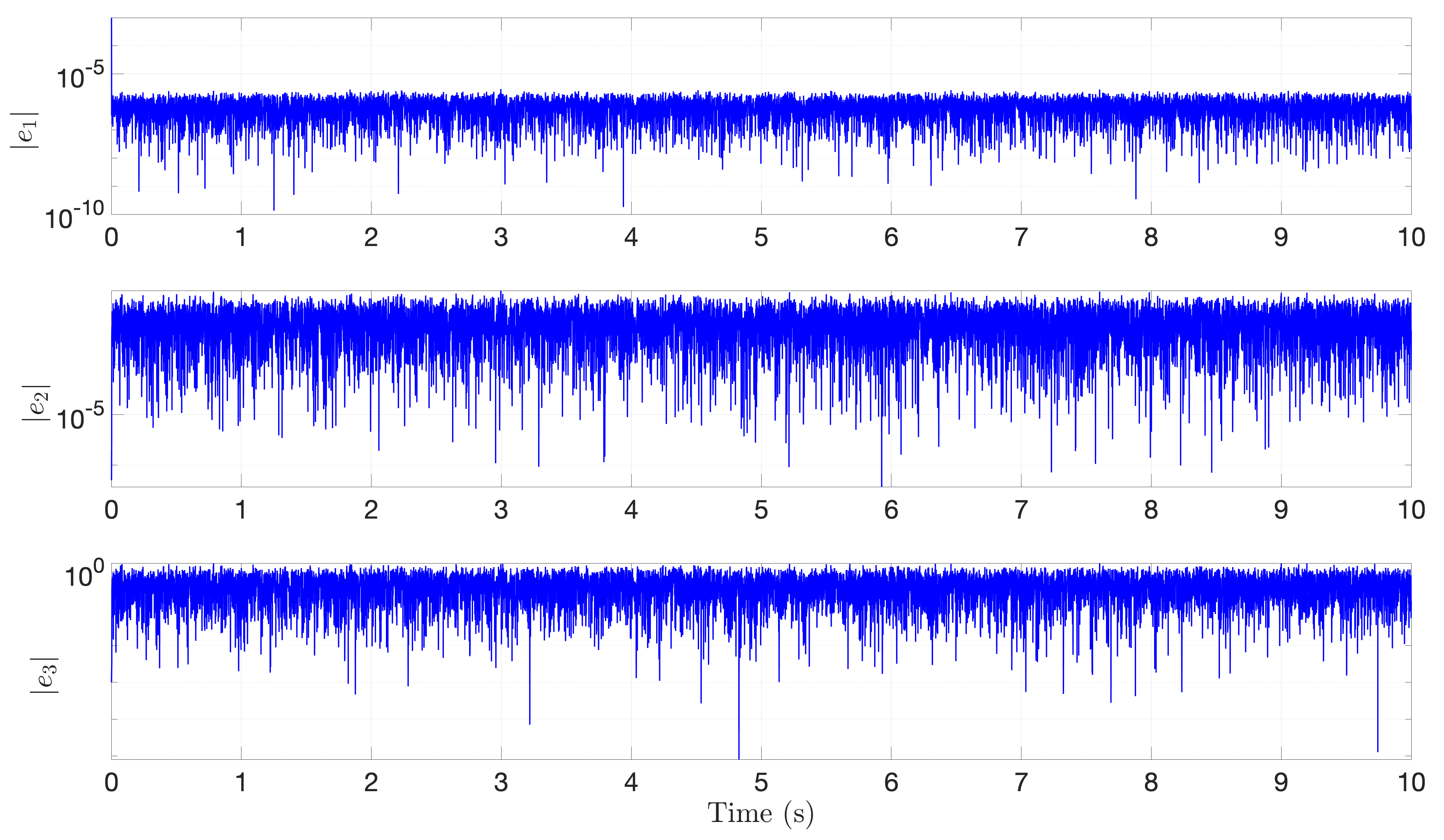}
      \end{center}
            \caption{Example of SIHD3-based differentiation of a sine function including extra noise - estimation error.}
       \label{fig_4}
\end{figure}

Low estimation errors are observed, thanks to the correction terms, without extra noise. In particular, the error is of order $< e_1 > = 2.96 \, 10^{-12}$  on the $z_1$ estimation, of order $< e_2 > = 4.69 \, 10^{-9}$ on the estimation of the first derivative $z_2$ and of order $< e_3 > = 4.49 \, 10^{-6}$ on the estimation of the second derivative $z_3$.
In presence of noise, the accuracy of the estimation remains correct for the two first orders ($< e_1 > = 7.55 \, 10^{-7}$ and $< e_2 > = 7.47 \, 10^{-4}$) but the noise becomes predominant for the third order, which makes the estimation very noisy. Note that the $\theta$ projector \cite{VSS_2022_michel}, which has shown interesting performances with respect to the rejection of the noise, has not been implemented in the current version of the proposed SIHD3 differentiation algorithm. It constitutes an important issue to deal with noisy signals for future investigations.

\clearpage
\section{Application to numerical integration scheme}

\subsection{Observer-based numerical scheme}

Consider an ordinary differential equation (ODE) such as:
\begin{equation}\label{eq:ODE_fond} 
\frac{d \, y(t)}{d \, t} = f(y(t), u(t)), \qquad t > t_0, \quad y(0) = y_0.
\end{equation}
\noindent
The quantities $u$ and $y$ represent respectively the input and the solution of (\ref{eq:ODE_fond}). The corresponding usual discrete explicit Euler scheme reads:

\begin{equation}\label{eq:EDP_discrete} 
\frac{y_{k+1} - y_k}{h} = f(y_k, u_k), \qquad k \in \mathbb{N}, \quad y(0) = y_0
\end{equation}

\noindent
where $k$ is the sampled time, $h$ is the time-step, and $y_k$, $y_{k+1}$ are respectively the solution of (\ref{eq:EDP_discrete}) at the discrete instants $t_k$ and $t_{k+1}$. The sampled are supposed equally distributed {\it i.e.} $t_k - t_{k-1}$ is constant for all $k$.

In the context of solving ODE, an observer-based numerical integration scheme is proposed to provide a better precision of the estimation of the ODE solution. The proposed scheme is close to the Taylor-based numerical schemes \cite{Griffiths, CASC_99} for which we point out the refined Taylor based estimation of $f(y_k, u_k)$ via our proposed SIHD algorithm. 

Our proposed differentiation method can be  associated to the nonstandard finite-difference (NSFD) methods \cite{Mickens_1, Mickens_2}, described as "powerful numerical methods that preserve significant properties of exact solutions of the corresponding differential equations" (a survey can be found {\it e.g.} in \cite{Patidar})\footnote{Applying a particular differentiator algorithm to the ODE discretization in the context of NFSD was proposed by the author in \cite{Michel2015}}.

The NSFD $\&$ SIHD3 based numerical integration scheme reads:

\begin{equation}\label{eq_nsfd_sihd}
\left\{\begin{array}{l}
z_1^+  = z_1 + \psi \, \frak{D}_{f}^{(1)}(z, z^+, \hat{y} ) \\
\hat{y} = z_1
\end{array}, \right. \, k \in \mathbb{N}, \, y(0) = y_0 \quad \mathrm{and} \quad \psi = h + O(h^2)
\end{equation}

\noindent
where the SIHD3 differentiator is used to estimate the $f$ part of the ODE and $\psi = h + O(h^2)$ that do satisfy $h \rightarrow 0$ according to the NSFD rules\footnote{For example, $\psi$ can be chosen linear of $h$ or exponential of $h$ \cite{Vigo}.}.

In this work, the purpose is to illustrate the convergence of the proposed NSFD $\&$  SIHD3-based algorithm by comparing the evolution of the error with the classical forward Euler scheme and the Runge-Kutta (RK) scheme according to the time.
Note that in this scheme, the semi-implicit differentiator  \eqref{eq:semi-implicit-D} appears as a closed-loop observer which tends to adjust and propagate the estimation of $f$ through the numerical integration scheme.

\subsection{Numerical examples and discussion}

Three examples are presented to illustrate some first convergence properties  through the evolution of the error $ | y_{exact}(t) - y_{ODE \, scheme}(t) | $ with the classical forward Euler scheme and the Runge-Kutta (RK) scheme according to the time. $y_{exact}(t)$ denotes the exact solution of the ODE \eqref{eq:ODE_fond} and $y_{ODE \, scheme}(t)$ denotes the approximated solution of the ODE \eqref{eq:EDP_discrete} using either the classical Euler forward method or the RK method.
 The internal parameters ($\lambda, \alpha$) of  SIHD3 are the same as in the first Section ($\lambda_1 = 10^3$, $\lambda_2 = 10^6$, $\lambda_3 = 10^9$, $\lambda_4 = 10^9$ and $\alpha = 0.95$).

\paragraph{Example 1} Consider the basic first order equation:

\begin{equation}\label{eq:ex_1} 
\frac{d \, y(t)}{d \, t} = -y(t) + 1, \quad y(0) = 10^{-2},  \quad t > 0
\end{equation}

\noindent
whose exact solution is $y_{exact} = 1 - \exp( - t )$. Fig. \ref{fig_nsfd_1} illustrates the comparaison of the errors between the schemes.
\begin{figure}[!b]
\begin{center}
   \includegraphics[width=12cm]{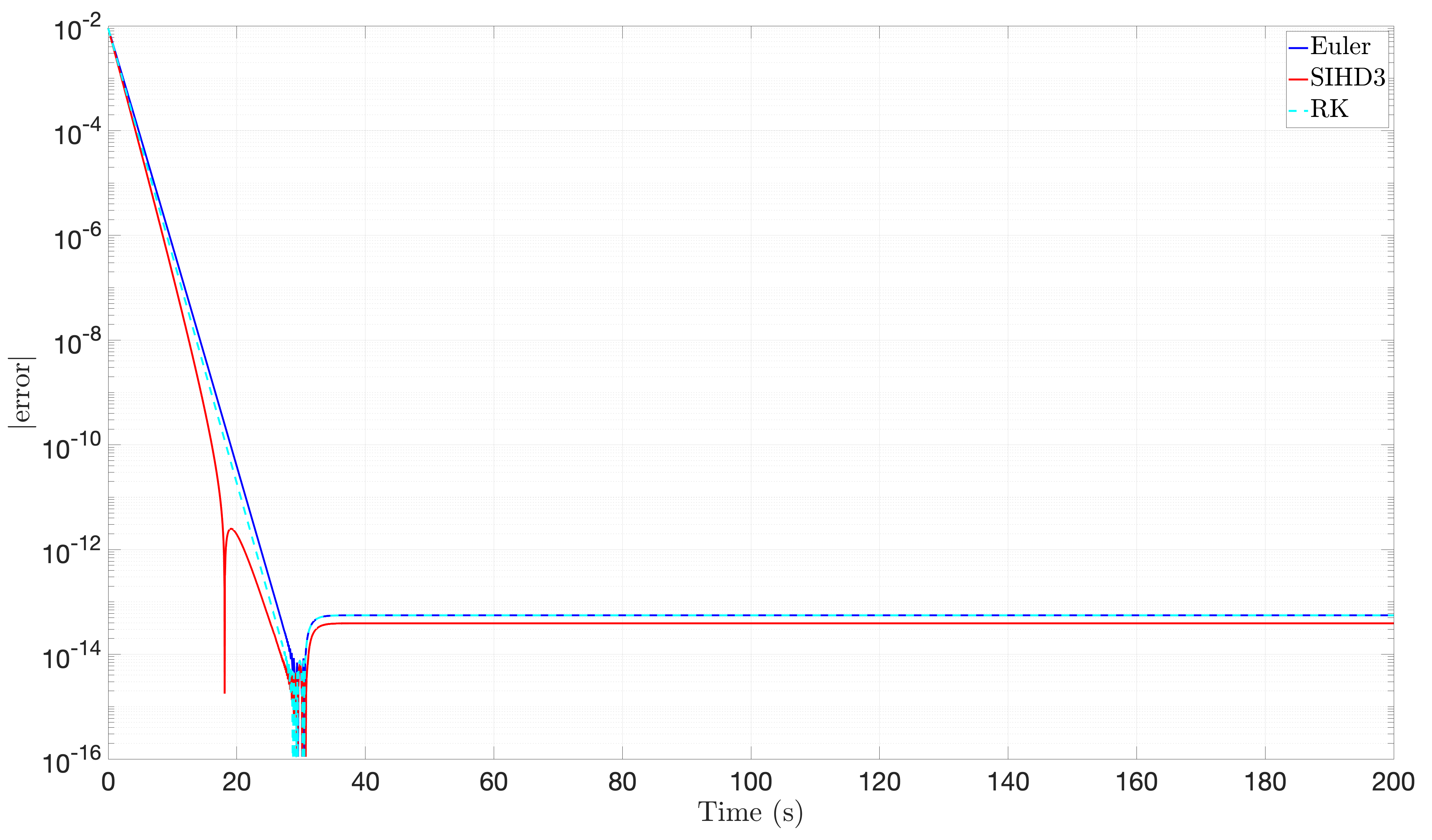}
      \end{center}
           \caption{Illustration of the SIHD3 scheme convergence compared with Euler forward scheme and Runge-Kutta scheme - Example 1.}
      \label{fig_nsfd_1}
\end{figure}

\paragraph{Example 2} Consider the logistic equation \cite{Griffiths}:

\begin{equation}\label{eq:ex_2} 
\frac{d \, y(t)}{d \, t} = 2 y ( 1 - y),  \quad y(0) = 1/5,  \quad t > 10
\end{equation}

\noindent
whose exact solution is $y_{exact} = \frac{1}{1 + 4 \exp( 2 (10-t  ))} $.  Fig. \ref{fig_nsfd_2} illustrates the comparaison of the errors between the schemes.
\begin{figure}[!b]
\begin{center}
   \includegraphics[width=12cm]{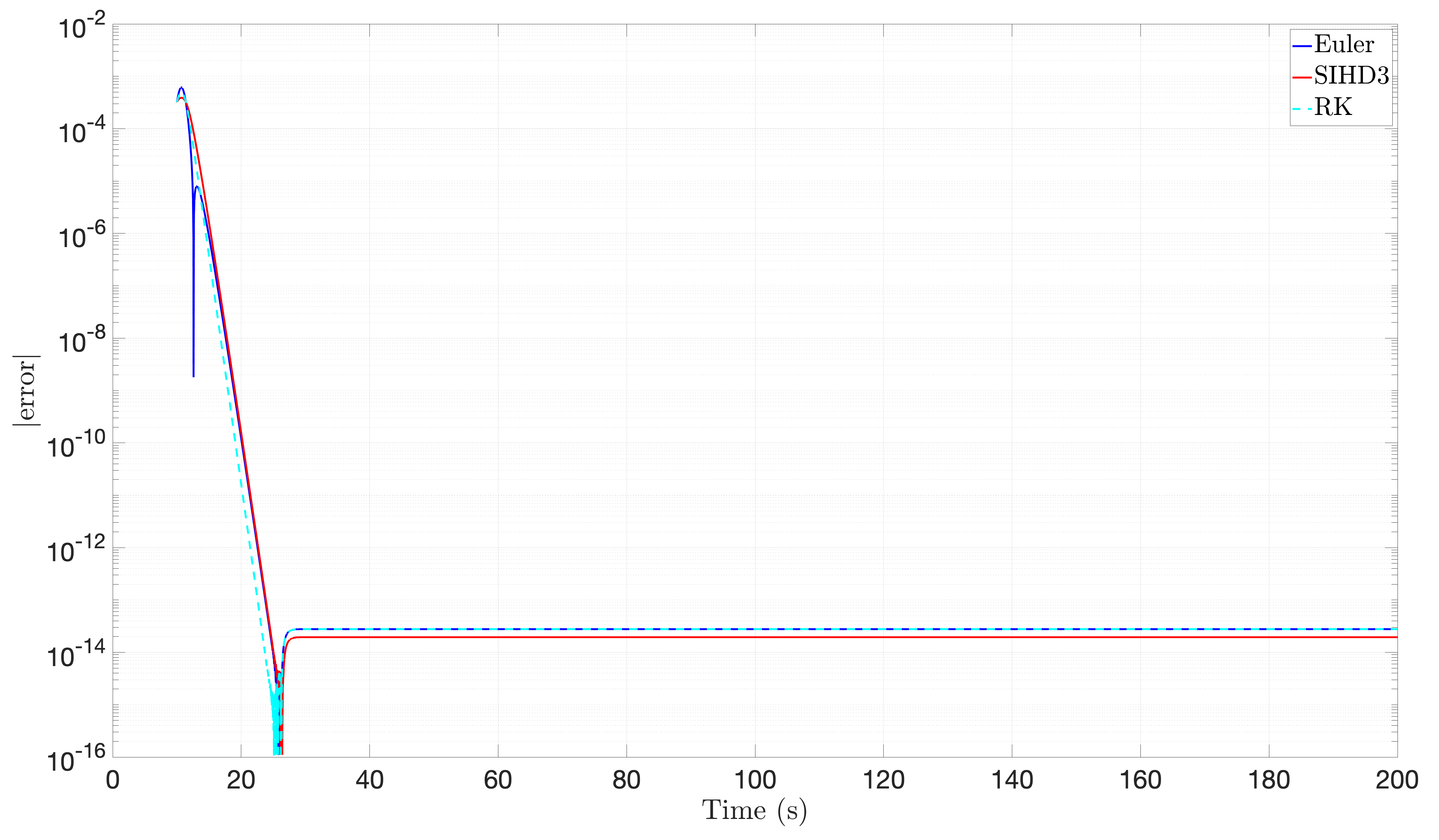}
      \end{center}
           \caption{Illustration of the SIHD3 scheme convergence compared with Euler forward scheme and Runge-Kutta scheme - Example 2.}
      \label{fig_nsfd_2}
\end{figure}

\paragraph{Example 3} Consider the simple model of a water-discharge problem \cite{book_Delft}:

\begin{equation}\label{eq:ex_3} 
\frac{d \, y(t)}{d \, t} = -10y^2 + 20,  \quad y(0) = 10^{-2},  \quad t > 0
\end{equation}

\noindent
whose exact solution is $y_{exact} = \sqrt{2} \tanh(\sqrt{300} t) $.  Fig. \ref{fig_nsfd_3} illustrates the comparaison of the errors between the schemes.
\begin{figure}[!b]
\begin{center}
   \includegraphics[width=12cm]{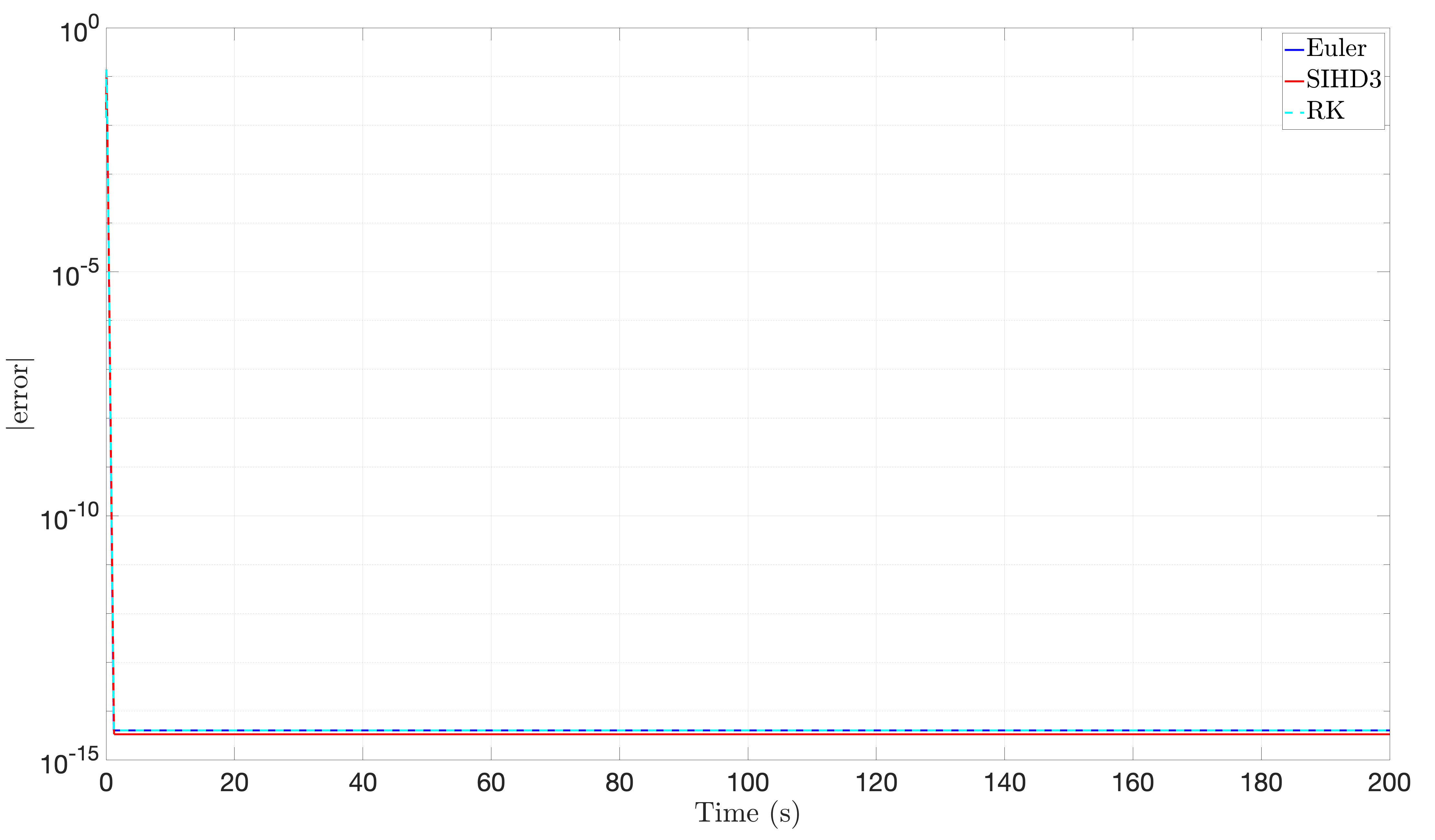}
      \end{center}
           \caption{Illustration of the SIHD3 scheme convergence compared with Euler forward scheme and Runge-Kutta scheme - Example 3.}
      \label{fig_nsfd_3}
\end{figure}

These preliminary numerical experiments show that the proposed NSFD-SIHD3 numerical scheme provides a slightly reduced error compared with the Euler and Runge-Kutta schemes, especially when the solution is converged towards the asymptotic solution. It is not easy to claim that the transient behavior is more precise than the other scheme since SIHD3 performs a dynamic estimation that depends on the ($\lambda, \alpha$) parameters, inducing delays during the estimation. Moreover, the NSFD approach appears naturally to adjust the rate of convergence of SIHD3 thanks to the $\psi$ time-step: indeed, $\psi$ can be considered as an additional SIHD3 'external' gain that can 'virtually' increase the speed of convergence of the scheme, hence giving some flexibility between the SIHD3 parameters ($\lambda, \alpha$) from the estimation side, and $\psi$ from the numerical scheme side.

\section{Conclusion and future works}

We have derived a   differentiator scheme of order three based on semi-implicit techniques that uses Taylor expansion in an implicit way. This differentiator allows reducing the estimation error in comparaison with exact derivatives. 
Future works include the implementation in real-time in the context of robotics applications.
We have also proposed a numerical scheme based on NSFD scheme including the  differentiator as a refined estimator of the $f$ term in the right side of the ODE.
Further investigations concern a complete study of convergence of the proposed scheme, as well as a comparison of the global error with existing algorithms and possibly an extension toward a semi-implicit scheme.

\section*{Appendix}

In this appendix, the corrections terms of the differentiator \eqref{eq:semi-implicit-diff-cor} are calculated from the assumption that the discretized signal $x_1$ can be differentiated until the third order (inducing a constant fourth order derivative).

Assuming that the discretized signal $x_1$ can be differentiated though the following 'multi-step' Taylor expansion:

\begin{eqnarray}
{x}_1^+&=& x_1+h\,x_2 + \frac{h^2}{2!} x_3 + \frac{h^3}{3!} x_4 + \frac{h^4}{4!} x_5 \nonumber \\
{x}_2^+&=& x_2+h\,x_3 + \frac{h^2}{2!} x_4 + \frac{h^3}{3!} x_5 \label{Syst_dis_Euler_y_exp}\\
{x}_3^+&=& x_3+h x_4 + \frac{h^2}{2!} x_5 \nonumber \\
{x}_4^+&=& x_4 + h x_5 \nonumber \\
{x}_5^+&=& x_5 \nonumber
\end{eqnarray}

\noindent
and assuming that the fifth order remains constant $x_5 = x_5^+$, the retro-propagation of the next-step predictions $x_q^+, q = {2...4}$ terms allows refining the first Taylor expansions by substituting each $x_q$ terms by their such corresponding predictions.

\noindent
To begin with, the third line:
\begin{eqnarray}
{x}_3^+ = x_3+h\,x_4 + \frac{h^2}{2!} x_5
\end{eqnarray}

\noindent
can be rewritten easily and substituting both $x_4 = x_4^+ - h x_5$ and $x_5^+ = x_5$, allows rewriting:
\begin{eqnarray*}
{x}_3^+ &=& x_3+h\, \textcolor{black}{( x_4^+ - h x_5^+ )} \textcolor{black}{ + \frac{h^2}{2!} x_5^+} 
\end{eqnarray*}

\noindent
that gives:
\begin{equation}
{x}_3^+ = x_3+h\, x_4^+ - \frac{h^2}{2!} x_5^+.
\end{equation}

In the same way, substituting $x_3 =  {x}_3^+ -h x_4^+ - \frac{h^2}{2!} x_5$ and $x_4 = x_4^+ - h x_5 $ in the second line gives after few calculations:


\begin{equation}
\label{x_2_correction_term}
x_2^+= x_2 + h \, x_3^+ \boldsymbol{- \frac{h^2}{2!} x_4^+} + R( x_5 )
\end{equation}

The term $R( x_5 )$ denotes the "residual" term associated to the useless fifth order.

Lastly, in the same way, substituting $x_3, x_4$ and  $x_2$ into $x_1$ gives: 

\begin{equation}
\label{x_1_correction_term}
{x}_1^+= x_1+h\, x_2^+ \boldsymbol{-  \frac{h^2}{2!} \, x_3^+ + \frac{h^3}{3!}  {x}_4^+  } + R( x_5^+ )
\end{equation}

The bold terms in \eqref{x_2_correction_term} and \eqref{x_1_correction_term} correspond indeed to the correction terms that refine the Taylor expansion of the two first lines in \eqref{eq:semi-implicit-diff-cor}, taking into account the refinement of the precision of the Taylor terms thanks to the "retro-propagation" of the high order implicit terms.

\section*{Acknowledgments}
This work was supported partially by the ANR project DigitSlid ANR-18-CE40-0008-01 (https://anr.fr/Project-ANR-18-CE40-0008). Jean-Pierre Barbot is supported with ({\it R\'egion Pays de la Loire}) Connect Talent GENYDROGENE project.

\bibliographystyle{unsrt}
\bibliography{BiblioJPB.bib, BiblioJit.bib, Bib_DIGITSLID_CDC.bib}

\end{document}